\documentclass[twoside,12pt]{amsart} 
\usepackage[dvips]{graphicx}
\DeclareGraphicsExtensions{ps,eps}
\makeatletter
\def\serieslogo@{}
\makeatother
\makeatletter 
\def\@setcopyright{}
\makeatother

\newtheorem{Theorem}{Theorem}[section]

\theoremstyle{definition}

\theoremstyle{remark}
\newtheorem{rem}{Remark}[section]

\newcommand{\rit}{{{\mathbb{R}}}}
\newcommand{\tit}{{{\mathbb{T}}}}
\newcommand{\zit}{{{\mathbb{Z}}}}

\begin{document}

\title
{
Relaxed solutions for incompressible inviscid flows: a variational and gravitational
approximation to the initial value problem
}

\author{Yann Brenier}
\address{
Y.B.: CNRS, D\'epartement de Math\'ematiques et Applications
\\Ecole Normale Sup\'erieure,  Universit\'e PSL, 45 rue d'Ulm 75005 Paris,
}
\author{Iv\'an Moyano}
\address{
I.M.: Laboratoire Jean Alexandre Dieudonn\'e
\\Universit\'e C\^ote-d'Azur,
Parc Valrose 06108 Nice Cedex 02.
}

\maketitle
\markboth{}{}

\section*{Abstract}
Following Arnold's geometric interpretation, the Euler equations of an incompressible fluid moving in a domain $D$ are known to be the optimality equation of the
minimizing geodesic problem along the group of orientation and volume preserving diffeomorphisms of $D$. This problem admits a well-established convex relaxation
which generates a set  of  ``relaxed", ``multi-stream",  version of the Euler equations. However, it is unclear that such relaxed equations are appropriate for the initial value problem
and the theory of turbulence, due to their lack of well-posedness for most initial data.  As an attempt to get a more relevant  set of relaxed Euler equations,
we address
 the multi-stream pressure-less gravitational Euler-Poisson system as an approximate model,
for which we show that the initial value problem can be stated as a concave maximization problem from which we can at least recover a large class
of smooth solutions for short enough times.
\section*{Introduction}
In \cite{YB-CPAM}, the first author introduced the following ``relaxed", ``multi-stream", version of the Euler equations for an incompressible
homogeneous fluid:
\begin{equation}\label{Euler-continuity}
\partial_t c +\nabla\cdot q=0,\;\;\;c=c(t,x,a)\ge 0,\;\;\;q=q(t,x,a)\in\rit^d,
\end{equation}
\begin{equation}\label{Euler-incompressibility}
\int_a c(t,x,a)=1,
\end{equation}
\begin{equation}\label{Euler-potential}
E=-\nabla\psi(t,x)\in\rit^d,\;\;\;
q=cv,\;\;\;v=v(t,x,a)=\nabla\theta(t,x,a)\in\rit^d,
\end{equation}
\begin{equation}\label{Euler-dynamic}
\partial_t (cv)+\nabla\cdot(c v\otimes v)=cE,
\end{equation}
where $t\in [0,T]$, $x\in D\subset \rit^d$, $a\in \mathcal A$ and $\nabla$ denotes the nabla operator
$$
\nabla=\left(\frac{\partial}{\partial{x^i}}\right)_{i=1,\cdot\cdot\cdot,d}\;\;{\rm{on}}\;\;\rit^d.
$$
The final time $T>0$ is fixed and, for simplicity, we assume the space domain $D$ to be the periodic cube $\tit^d=(\rit/\zit)^d$.
The space of labels $\mathcal A$ is taken as a compact metric space equipped with a Borel probability measure $\mu$
(typically, the counting measure when $\mathcal A$ is discrete and finite, or the Lebesgue measure when $\mathcal A=D$)
and  the integration in $a\in\mathcal A$ in formula (\ref{Euler-incompressibility}) is performed according to $\mu$.
This corresponds to the description of an incompressible fluid viewed as a superposition of
different streams, labelled by $a$, moving across each other with their own concentration $c=c(t,x,a)\ge 0$, potential velocity $v=v(t,x,a)=\nabla\theta(t,x,a)\in \rit^d$
and momentum $q(t,x,a)=c(t,x,a)v(t,x,a)\in\rit^d$
fields, while they are driven by a $common$ potential acceleration field
$E=-\nabla p(t,x)$, independent of $a$, which maintains the equal
occupation of volumes through (\ref{Euler-incompressibility}).
When $\mathcal A$ is taken to be $D$ itself, $\mu$ being the Lebesgue measure, these equations turn out to to be the correct ``relaxation" of the Euler equations when solving the minimizing geodesic problem (MGP)
in the framework of their geometric interpretation, going back to Arnold \cite{Arnold},
in terms of the geodesic flow over the group $SDiff(D)$  of all orientation and volume-preserving diffeomorphisms of $D$. Let us give a short explanation of this statement.
The MGP can definitely be uniquely solved, in the classical setting of diffeomorphisms, in a small Sobolev neighborhood of the identity map, as shown by
Ebin and Marsden \cite{EM}, but not in the large as shown by Shnirelman \cite{Shni}, at least as $d\ge 3$. Indeed, minimizing sequences may develop small scale oscillations and
admit no limit in the classical setting. Therefore, one has to complete $SDiff(D)$ by the set of ``doubly stochastic measures" $DS(D)$ (also called ``polymorphisms" as in
\cite{Neretin}), i.e. the set all Borel measures $c$ over $D\times D$ that project to the Lebesgue measure on each copy of $D$,
each diffeomorphism
\\
$a\in D\rightarrow X(a)\in D$ of $SDiff(D)$ generating a corresponding $c$ through 
\\
$c(x,a)=\delta(x-X(a))$ or, more precisely, 
$$
\int_{D\times D}f(x,a)c(dxda)=\int_D f(X(a),a)da,\;\;\;\forall f\in C(D\times D).
$$
This completion process can be justified in a static way, as in \cite{Neretin}, or, even better, in a dynamical way as in \cite{Shni-GAFA}.
Then, as shown in \cite{YB-CPAM,AF1}, the MGP can be ``relaxed" as a convex minimization problem over $DS(D)$.
The optimality equations turn out to be the ``multi-stream" equations
(\ref{Euler-continuity},\ref{Euler-incompressibility},\ref{Euler-potential},\ref{Euler-dynamic}), 
in the special case when $\mathcal A$ is taken as $D$.
Nevertheless, the ``relaxed" MGP, for which time-boundary data $c(0,x,a)$ and $c(T,x,a)$ are provided at both $t=0$ and $t=T$, without any information on the velocity
field, differs very much from the initial value problem (IVP), when both $c(0,x,a)$ and $v(0,x,a)$ are prescribed at the initial time $t=0$, without any information on the final time $T$.
Therefore, it is not clear that the ``relaxed Euler equations" 
(\ref{Euler-continuity},\ref{Euler-incompressibility},\ref{Euler-potential},\ref{Euler-dynamic})
are relevant as a good relaxation of
the Euler equations for the IVP.
Notice they can also be written in ``Vlasov form" 
\begin{equation}\label{VE}
\partial_t f(t,x,\xi)+\nabla_x\cdot (\xi f(t,x,\xi))+\nabla_\xi\cdot (E(t,x)f(t,x,\xi))=0,\;\;(t,x,\xi)\in [0,T]\times D\times \rit^d,
\end{equation}
\begin{equation}\label{VE-Poisson}
\int_{\xi\in\rit^d} f(t,x,\xi)=1,\;\;\;
E=-\nabla\psi(t,x)\in\rit^d,
\end{equation}
where
\begin{equation}\label{MEPS-VP}
f(t,x,\xi)=\int_{a\in\mathcal A}c(t,x,a)\delta(\xi-v(t,x,a)).
\end{equation}
These equations can be seen as a ``kinetic formulation" of the Euler equations. Observe that, after integration in $\xi\in\rit^d$, we get
for the ``macroscopic" quantities
$$
V(t,x)=\int_{\xi\in\rit^d} \xi f(t,x,\xi),\;\;\;M(t,x)=\int_{\xi\in\rit^d} \xi\otimes\xi f(t,x,\xi),
$$
the (incomplete) set of equations
$$
\nabla\cdot V=0,\;\;\;\partial_t V+\nabla\cdot M+\nabla\psi=0,
$$
together with inequality $M\ge V\otimes V$, in the sense of symmetric matrices. In other words, $(V,M)$ is a $subsolution$ to the Euler
equations, in the De Lellis-Sz\'ekelyhidi framework of ``convex integration" \cite{DS}. It is known that such subsolutions, at least when they are strict, i.e. $M>V\otimes V$, can be approximated by standard weak solutions to the Euler equations thanks to ``convex integration" tools \cite{DS}.
(We also refer to \cite{YB-CPAM,YB-DS}
for a discussion about the connections between such a Vlasov-type formulation, the concept of ``sharp measure-valued solutions"
to the Euler equations and the older concept of ``measure-valued solutions".)
\\
It is important to notice that equations (\ref{VE},\ref{VE-Poisson}) can also be thought as the highly singular limit of
the well-known gravitational Vlasov-Poisson system,
\begin{equation}\label{VP}
\partial_t f(t,x,\xi)+\nabla_x\cdot (\xi f(t,x,\xi))+\nabla_\xi\cdot (E(t,x)f(t,x,\xi))=0,\;\;(t,x,\xi)\in [0,T]\times D\times \rit^d,
\end{equation}
\begin{equation}\label{VP-Poisson}
1-\epsilon \nabla\cdot E(t,x)=
\int_{\xi\in\rit^d} f(t,x,\xi),\;\;\;
E=-\nabla\psi(t,x)\in\rit^d,
\end{equation}
as $\epsilon\downarrow 0$.
In sharp contrast with the Vlasov-Poisson system, for which the IVP is well-posed,
the limit equations are not well-posed in the usual sense \cite{Han-Kwan}, although there is a set of initial conditions, defined by a suitable ``Penrose-Rayleigh" type
condition around which the IVP is presumably well-posed \cite{YB-NL}. As a matter of fact, the study of the limiting process is very delicate and is an
active field of research, for which we refer, as a very recent reference, to \cite{Baradat}.
The Vlasov-Poisson system itself admits a ``multi-stream formulation" which has been used for a while, for instance by Grenier in \cite{Grenier} and recently
by Baradat \cite{Baradat} (at least in the framework of Plasma Physics for which $\epsilon$ is negative). The corresponding
``multi-stream" (pressure-less) Euler-Poisson system (MEPS) reads
\begin{equation}\label{MEPS-continuity}
\partial_t c +\nabla\cdot q=0,\;\;\;c=c(t,x,a)\ge 0,\;\;\;q=q(t,x,a)\in\rit^d,
\end{equation}
\begin{equation}\label{MEPS-Poisson}
1-\epsilon \nabla\cdot E(t,x)=
\int_a c(t,x,a),
\end{equation}
\begin{equation}\label{MEPS-potential}
E=-\nabla\psi(t,x)\in\rit^d,\;\;\;
q=cv,\;\;\;v=v(t,x,a)=\nabla\theta(t,x,a)\in\rit^d,
\end{equation}
\begin{equation}\label{MEPS-dynamic}
\partial_t (cv)+\nabla\cdot(c v\otimes v)=cE.
\end{equation}
This system is (formally) energy conservative:
\begin{equation}\label{MEPS-energy}
\frac{d}{dt}\left(\int_{(x,a)\in D\times \mathcal A}c(t,x,a)\frac{|v(t,x,a)|^2}{2}-\int_{x\in D}\epsilon\frac{|\nabla\psi(t,x)|^2}{2}\right)=0.
\end{equation}
Let us now sketch a physical interpretation of these equations.
The MEPS describes the dynamics of a self-gravitating pressure-less multi-stream fluid,  moving according to the standard Newtonian model of gravitation.  
\\
Each stream, labelled by some $a$ in $\mathcal A$,
admits its own concentration field $c(t,x,a)\ge 0$ and potential velocity field $v(t,x,a)=\nabla\theta(t,x,a)\in\rit^d$, while the momentum $c(t,x,a)v(t,x,a)$
is denoted by $q(t,x,a)$. All together,  these streams share a common
acceleration field $E(t,x)=-\nabla\psi(t,x)$ where the gravitational potential $\psi$ is coupled to them through the Poisson equation
\begin{equation}\label{Poisson}
\-\epsilon \Delta\psi(t,x)=1-\int_a c(t,x,a),
\end{equation}
with $\epsilon>0$, $\epsilon^{-1}$ being the gravitational constant, properly rescaled. We assume the initial concentration field $c_0$ 
to satisfy
$$
\int_{(x,a)\in  D\times \mathcal A} c_0(x,a)=1,
$$
so that the ``total mass" is always equal to 1:
\begin{equation}\label{MEPS-mass}
\int_{(x,a)\in  D\times \mathcal A} c(t,x,a)=1,
\end{equation}
because of (\ref{MEPS-continuity}) (since $D$ is periodic). 
This way,  in the Poisson equation (\ref{Poisson}), the right-hand side has zero mean over $D$, where
constant 1 induces  a background repulsive gravitational 
potential due to the periodicity of the spatial domain $D$ (as usual in computational Cosmology \cite{Frisch2,YB-Frisch}).
$$
$$
Ultimately, our hope would be to adapt the relaxation technique successfully used for the MGP to address the IVP in the case of the Euler equations, with the goal of getting
new ``relaxed" equations, that might be of some interest for the 
theory of turbulence \cite{Frisch}, beyond the one we started with, namely 
(\ref{Euler-continuity},\ref{Euler-incompressibility},\ref{Euler-potential},\ref{Euler-dynamic}).
As a preliminary step, we address in
the present paper the case of the more accessible, but nevertheless interesting, gravitational Vlasov-Poisson system, with finite $\epsilon>0$. In that case,
we first explain, in the spirit of \cite{YB-CPAM,Loeper},  how the solutions of the MEPS 
can be recovered on a given time interval $[0,T]$ by the least action principle through a suitable space-time convex minimization problem, the concentration fields being prescribed at both $t=0$ and $t=T$ without any information required on the velocity field. 
Next, we show that, surprisingly enough, through a suitable $augmented$  $Lagrangian$ technique, the initial value problem, when both concentration and velocity fields
are prescribed at $t=0$, without any information needed at time $t=T$,
can $also$ be relaxed as a (dual) concave maximization problem,  through Theorem \ref{dual}. We finally explain how, under some smallness
condition on $T$, smooth solutions of the MEPS system can be recovered, through Theorem \ref{nogap}.
Let us conclude this introduction by saying that the results obtained in the present paper for the multi-stream Euler-Poisson system MEPS 
for each fixed $\epsilon>0$ give some hope that a similar analysis could be performed in the
more challenging case of the multi-stream Euler equations as $\epsilon\downarrow 0$.
\section{The Least Action Principle for the MEPS}
As already established by Loeper in \cite{Loeper} in the case of a single stream, the MEPS can be derived from
the $convex$ minimization problem in $(c, q)$
\begin{equation}\label{MEPS-LAP}
\inf _{c,q,E}\;\int_{Q'} 
\frac{|q|^2}{2c}+\int_{Q}\frac{\epsilon|E|^2}{2},
\end{equation}
where
\begin{equation}\label{Q'}
Q=[0,T]\times D,
\;\;\;Q'=Q\times \mathcal A,
\end{equation}
as the fields $c$, $q$, $E$ are subject to $linear$ constraints
(\ref{MEPS-continuity},\ref{MEPS-Poisson})
and $c$ is prescribed at both $t=0$ and $t=T$.
Indeed, using Lagrange multipliers for constraints (\ref{MEPS-continuity},\ref{MEPS-Poisson}),
we first get the equivalent saddle-point formulation of the minimization problem:
\begin{equation}\label{MEPS-saddle}
\inf _{c,q,E}\;\sup _{\theta,\psi}\;
BT_0(\theta)+\int_{(t,x)\in Q} \psi+\epsilon E\cdot \nabla\psi+\frac{\epsilon|E|^2}{2}
+\int_{(t,x,a)\in Q'} \frac{|q|^2}{2c}
-(\partial_t \theta+\psi)c-q\cdot\nabla\theta,
\end{equation}
where
\begin{equation}\label{MEPS-BT}
BT_0(\theta)=\int_{(x,a)\in D\times\mathcal A}c_T(x,a)\theta(T,x,a)-c_0(x,a)\theta(0,x,a).
\end{equation}
Then, differentiating the Lagrangian with respect to $c$, $q$ and $E$, we get
$$
q=c\nabla\theta,\;\;\;-\frac{|q|^2}{2c^2}-\partial_t\theta-\psi=0,\;\;\;E=-\nabla\psi,
$$
which implies
\begin{equation}\label{MEPS}
\partial_t\theta(t,x,a)+\frac{|\nabla\theta(t,x,a)|^2}{2}+\psi(t,x)=0,\;\;\;\epsilon\Delta\psi(t,x)=\int_a c(t,x,a)-1,
\end{equation}
leading to (\ref{MEPS-potential}) and (\ref{MEPS-dynamic})
using (\ref{MEPS-continuity}).

\section{Saddle-point formulation of the initial value problem (IVP)}
We now address the IVP by
$adding$ to the Lagrangian the extra term
$$
\int_{(t,x,a)\in Q'}
-\partial_t A\cdot q
-\nabla A\cdot(\frac{q\otimes q}{c})-cA\cdot E
-\int_{(x,a)\in D\times\mathcal A}A(0,x,a)\cdot q_0(x,a),
$$
which takes into account the $weak$ formulation of  (\ref{MEPS-dynamic}), where
$A=A(t,x,a)\in\rit^d$ is any test function subject to $A(T,x,a)=0$.
Then, we get a new saddle-point
problem
\begin{equation}\label{IVP-saddle}
I(c_0,q_0)=\inf _{c,q,E}\;\sup _{\theta,A,\psi}\;
BT(\theta,A)+\int_{(t,x)\in Q} \psi+\epsilon E\cdot \nabla\psi+\frac{\epsilon|E|^2}{2}
\end{equation}
$$
+\int_{(t,x,a)\in Q'} \frac{|q|^2}{2c}
-(\partial_t \theta+\psi)c-q\cdot\nabla\theta
-\partial_t A\cdot q-\nabla A\cdot(\frac{q\otimes q}{c})-cA\cdot E
$$
where
\begin{equation}\label{IVP-BT}
BT(\theta,A)=-\int_{(x,a)\in D\times\mathcal A}c_0(x,a)\theta(0,x,a)+q_0(x,a)\cdot A(0,x,a)
\end{equation}
and test functions $\theta$ and $A$ must vanish at time $t=T$.
This ``augmented Lagrangian" strategy allows us to input initial condition $q_0$, together with $c_0$,
while data $c_T$ is no longer needed. This idea has been already applied by the first author to the Euler equations
of incompressible fluids and, also, to the class of first-order systems of conservation laws with 
a convex entropy in \cite{YB-CMP} and extended to various interesting models (including ideal MHD) 
by Vorotnikov in \cite{Vorot}.

\section{Dual formulation of the IVP for the MEPS}
By exchanging the sup and the inf, we obtain the "dual" problem 
\begin{equation}\label{IVP-dual0}
J(c_0,q_0)
=\sup _{\theta,A,\psi}\;\inf _{c,q,E}\;(\cdot).
\end{equation}
Observe that we have a priori no more than the ``weak" duality property
$$
I(c_0,q_0)=\inf _{c,q,E}\;\sup _{\theta,A,\psi}\;(\cdot)\;\;\ge\;\sup _{\theta,A,\psi}\;\inf _{c,q,E}\;(\cdot)\;=J(c_0,q_0).
$$
Indeed, by augmenting the Lagrangian, we have destroyed the convex structure of the original problem,
because of the nonlinear term $q\otimes q/c$ and, therefore, a duality gap cannot be excluded.
Anyway, the dual problem corresponds to
a $concave$ maximization problem in $(\theta,A,\psi)$, due to the linearity of
weak formulations with respect to their test functions. Let us compute it more explicitly.
First, we perform the infimum in $q$. We see that this infimum is $-\infty$ unless 
$A$ satisfies the inequality constraint
\begin{equation}\label{IVP-inequality1}
(\nabla A+\nabla A^T)(t,x,a)\le I,
\end{equation}
pointwise in the sense of symmetric matrices, where $I$ denotes the identity $d\times d$ matrix.
In that case, the optimal value of $q$ is obtained as
$$
(I-\nabla A-\nabla A^T)q=(\partial_t A+\nabla\theta)c
$$
and we find, after minimization in $q$,
\begin{equation}\label{IVP-dualbis}
J(c_0,q_0)=\sup _{\theta,A,\psi}
BT(\theta,A)+\;\inf _{c\ge 0,E}
\int_{Q} \psi+\epsilon E\cdot \nabla\psi+\frac{\epsilon|E|^2}{2}
\end{equation}
$$
-\int_{Q'} (\frac{1}{2}(I-\nabla A-\nabla A^T)^{-1}(\partial_t A+\nabla\theta)^{\otimes 2}+\partial_t\theta+\psi+E\cdot A)c.
$$
(Here notation $Mw^{\otimes 2}$ stands for $w\cdot(Mw)$ whenever $M$ is a
$d\times d$ matrix and $w$ is a vector in $\rit^d$.) Let us now perform the infimum in $c\ge 0$. We first observe that
$$
\inf _{c\ge 0,E}\le \inf _{c\ge 0,E=0}
$$
which is $-\infty$ unless
\begin{equation}\label{IVP-inequality2}
-\eta(t,x,a)=
(\frac{1}{2}(I-\nabla A-\nabla A^T)^{-1}(\partial_t A+\nabla\theta)^{\otimes 2}+\partial_t\theta)(t,x,a)+\psi(t,x)\le 0.
\end{equation}
Thus, we may now perform the infimum in $c$, under this second inequality constraint, and immediately obtain
\begin{equation}\label{IVP-dualter}
J(c_0,q_0)=\sup _{\theta,A,\psi}
BT(\theta,A)+\;\inf _{E}
\int_{Q} \psi+\frac{\epsilon|E|^2}{2}+\epsilon E\cdot \nabla\psi
\end{equation}
where $E$ is subject to the pointwise inequality
$$
E(t,x)\cdot A(t,x,a)\le \eta(t,x,a).
$$
So, we have finally obtained:
\begin{Theorem} 
\label{dual}
The concave dual maximization problem, proposed to solve the IVP for the MEPS
(\ref{MEPS-continuity},\ref{MEPS-Poisson},\ref{MEPS-potential},\ref{MEPS-dynamic}) with initial condition $(c_0,q_0)$, 
reads:
\begin{equation}\label{IVP-dual}
J(c_0,q_0)=\sup _{\theta,A,\psi}
BT(\theta,A)+\int_{(t,x)\in Q} \left(\psi(t,x)-\epsilon K_{(A,\eta)(t,x,\cdot)}(\nabla\psi(t,x))\right)dxdt,
\end{equation}
where $(\theta,A)$ must vanish at $t=T$, $A$ is subject to
$(\nabla A+\nabla A^T)(t,x,a)\le I$, 
\begin{equation}\label{IVP-eta}
 \eta(t,x,a)=-(\frac{1}{2}(I-\nabla A-\nabla A^T)^{-1}(\partial_t A+\nabla\theta)^{\otimes 2}+\partial_t\theta)(t,x,a)-\psi(t,x)
\end{equation} 
must be nonnegative, and we use notations:
$$
BT(\theta,A)=-\int_{(x,a)\in D\times\mathcal A}c_0(x,a)\theta(0,x,a)+q_0(x,a)\cdot A(0,x,a),
$$
$$
K_{(A,\eta)(\cdot,\cdot,\cdot)}(B)=
-\inf\{\;\frac{|E|^2}{2}+ E\cdot B,\;\;E\in\rit^d\;s.\;t.\;\;E.A(\cdot,\cdot,a)\le \eta(\cdot,\cdot,a),\;\;\forall a\in\mathcal A\}.
$$
\end{Theorem}
\section{A ``no duality-gap" result}
\begin{Theorem} 
\label{nogap}
Let $(c^s>0,q^s,E^s)$ be a smooth solution  on $Q'=[0,T]\times D\times\mathcal A$ to the MEPS
of form
$$
q^s=c^sv^s,\;\;\;v^s(t,x,a)=\nabla\theta^s(t,x,a),\;\;\;E^s(t,x)=-\nabla\psi^s(t,x).
$$
We make two assumptions on such a solution. First,
\begin{equation}\label{IVP-Ponce}
(\nabla v^s+(\nabla v^s)^T)(t,x,a)< \frac{I}{T-t},
\end{equation}
holds true pointwise in $Q'$, in the sense of symmetric matrices, where $I$ denotes the identity $d\times d$ matrix.
Second, the velocity field $v^s$ is ``weakly absorbing" in the sense that, for all fixed $(t,x)\in Q$ and for all vector
$V$ in $\rit^d$, there is a nonnegative mesure $\lambda$ (depending on $t,x$ and $V$) on $\mathcal A$ such that
\begin{equation}\label{IVP-absorbing}
V=\int_{a\in\mathcal A} v^s(t,x,a)\lambda(da).
\end{equation}
Then there is no duality gap between the primal and dual problems, namely
$$
I(c_0,q_0)=J(c_0,q_0).
$$
Moreover, a solution of the dual problem (\ref{IVP-dual}) is explicitly given by
\begin{equation}\label{IVP-dual-solution}
A(t,x,a)=(t-T)v^s(t,x,a), \;\;\;\psi(t,x)=\partial_t((T-t)\psi^s(t,x)),
\end{equation}
$$
\;\;\; \theta(t,x,a)=(t-T)(\psi^s(t,x)-\frac{|v^s(t,x,a)|^2}{2}).
$$
\end{Theorem}
N.B. In these notations, the superscript ``$s$" stands for ``solution" and enables us to make a clear and crucial
distinction between the solution and the various test functions used in the following calculations.
\begin{rem}
Notice that assumption (\ref{IVP-absorbing}) is automatically satisfied as soon as, for each $(t,x)$ the convex hull of the velocity 
range $\{v^s(t,x,a),\:\:a\in\mathcal A\}$ contains a small open ball containing $0$ in $\rit^d$. Surprisingly enough, in the case 
of a single stream, i.e. when there is only one label $a$,  condition (\ref{IVP-absorbing}) implies $v^s=0$
and rules out all non trivial solutions!
\end{rem}
\subsection*{Proof of Theorem \ref{nogap}}
\subsubsection*{Step 1}

Since $(c^s,q^s,E^s)$ is a smooth solution to the MEPS with initial condition $(c_0,q_0)$, it weakly solves both (\ref{MEPS-continuity}) and (\ref{MEPS-dynamic}).
Thus, for all $(\theta,A)$ that vanish at $t=T$, we get
$$
BT(\theta,A)
+\int_{Q'}-\partial_t \theta c^s-q^s\cdot\nabla\theta
-\partial_t A\cdot q^s-\nabla A\cdot(\frac{q^s\otimes q^s}{c^s})-c^sA\cdot E^s=0,
$$
where we recall that the time-boundary term $BT$ is defined by (\ref{IVP-BT}).
We also have (\ref{MEPS-Poisson}) which implies, for all $\psi$:
$$
\int_{(t,x)\in Q} \psi(t,x)+\epsilon E^s\cdot \nabla\psi(t,x)
-\int_{(t,x,a)\in Q'}\psi(t,x) c^s(t,x,a)=0.
$$
So, the supremum in $(\theta,A,\psi)$ of 
$$
BT(\theta,A)+\int_{Q} \psi^s+\epsilon E^s\cdot \nabla\psi^s+\frac{\epsilon|E^s|^2}{2}
$$
$$
+\int_{Q'} \frac{|q^s|^2}{2c^s}
-(\partial_t \theta+\psi)c^s-q^s\cdot\nabla\theta
-\partial_t A\cdot q^s-\nabla A\cdot(\frac{q^s\otimes q^s}{c^s})-c^sA\cdot E^s,
$$
is just equal to
$$
\int_Q\frac{\epsilon|E^s|^2}{2}+\int_{Q'} \frac{|q^s|^2}{2c^s},
$$
where we recall that $q^s=cv^s$, $v^s=\nabla\theta^s$ and $E^s=-\nabla\psi^s$.
By definition (\ref{IVP-saddle}) of the primal problem, this supremum
is an upper bound for $I(c_0,q_0)$ and, therefore, for $J(c_0,q_0)$ as well.
Thus, at the end of this first step, we already have
\begin{equation}\label{upper}
\int_Q\frac{\epsilon|\nabla\psi^s|^2}{2}+\int_{Q'} \frac{c|v^s|^2}{2}\ge I(c_0,q_0)\ge J(c_0,q_0). 
\end{equation}
\subsubsection*{Step 2}

Let us now move to the dual side.
By definition (\ref{IVP-dual}), we get the lower bound
$$
J(c_0,q_0)\ge 
BT(\theta,A)+
\int_{Q} \psi
+\;\;\inf _{c,q,E}\;\;
\int_{Q} \epsilon E\cdot \nabla\psi+\frac{\epsilon|E|^2}{2}
$$
$$
+\int_{Q'} \frac{|q|^2}{2c}
-(\partial_t \theta+\psi)c-q\cdot\nabla\theta
-\partial_t A\cdot q-\nabla A\cdot(\frac{q\otimes q}{c})-cA\cdot E,
$$
whatever is our choice of $(\theta,A,\psi)$.
Let us make the ``ansatz"
\begin{equation}\label{ansatz}
A=(t-T)v^s(t,x), \;\theta=(t-T)(\psi^s(t,x)-\frac{|v^s(t,x,a)|^2}{2}),
\end{equation}
while, the choice of $\psi$ will be made later.
So, $(A,\theta,\psi)$ being fixed, we may write our lower bound as
\begin{equation}\label{lower}
J(c_0,q_0)\ge \mathcal J+
BT(\theta,A)+\int_{Q} \psi,
\end{equation}
where we denote
\begin{equation}\label{J}
\mathcal J=
\inf _{c,q,E}\int_{Q} \epsilon E\cdot \nabla\psi+\frac{\epsilon|E|^2}{2}
\end{equation}
$$
+\int_{Q'} \frac{|q|^2}{2c}
-(\partial_t \theta+\psi)c-q\cdot\nabla\theta
-\partial_t A\cdot q-\nabla A\cdot(\frac{q\otimes q}{c})-cA\cdot E.
$$
Let us first focus on the infimum in $q$ in this definition of $\mathcal J$.
Since we have assumed $(T-t)(\nabla v^s+(\nabla v^s)^T)< I$,  inequality $(I-\nabla A-\nabla A^T)>0$ holds true in the sense of symmetric
matrices and, therefore, the infimum in $q$ is obtained for
$$
(I-\nabla A-\nabla A^T)q=(\partial_t A+\nabla\theta)c
$$
i.e., setting $q=cv$, as
$$
(I-\nabla A-\nabla A^T)v=\partial_t A+\nabla\theta,
$$
which reads, by definition of $(\theta,A)$:
$$
(I-(t-T)\nabla v^s-(t-T)(\nabla v^s)^T)v=v^s+(t-T)\partial_t v^s+(t-T)\nabla(\psi^s-\frac{|v^s|^2}{2}),
$$
or, equivalently,
$$
v-v^s=(t-T)\left(\partial_t v^s+(\nabla v^s+(\nabla v^s)^T)v+\nabla(\psi^s-\frac{|v^s|^2}{2})\right).
$$
Observe that this equation admits $v=v^s=\nabla\theta^s$ as an obvious solution. Indeed, we have
$$
\partial_t v^s+(\nabla v^s+(\nabla v^s)^T)v^s+\nabla(\psi^s-\frac{|v^s|^2}{2})
=\nabla\left(\partial_t\theta^s+\frac{|\nabla\theta^s|^2}{2}+\psi^s\right)=0
$$
since $\theta^s$ and $\psi^s$ solve (\ref{MEPS}).
At this stage, we may therefore rewrite $\mathcal J$ defined by (\ref{J}) as
$$
\mathcal J=
\inf _{c\ge 0,E}\int_{Q} \epsilon E\cdot \nabla\psi+\frac{\epsilon|E|^2}{2}
$$
$$
+\int_{Q'}\left(\frac{|v^s|^2}{2}
-\partial_t \theta-\psi-v^s\cdot\nabla\theta
-\partial_t A\cdot v^s-\nabla A\cdot(v^s \otimes v^s)-A\cdot E\right)c.
$$
The infimum with respect to $c\ge 0$ immediately leads to the new value
of $\mathcal J$:
\begin{equation}\label{Jbis}
\mathcal J=
\inf\int_{Q} \epsilon E\cdot \nabla\psi+\frac{\epsilon|E|^2}{2}
\end{equation}
where $E$ is subject to the pointwise inequality
$$
\frac{|v^s|^2}{2}
-\partial_t \theta-\psi-v^s\cdot\nabla\theta
-\partial_t A\cdot v^s-\nabla A\cdot(v^s \otimes v^s)-A\cdot E\ge 0.
$$
Let us recall that $\theta,A$ have been fixed
according to  (\ref{ansatz}), i.e. :
$$
A=(t-T)v^s, \;\;\theta=(t-T)(\psi^s-\frac{|v^s|^2}{2}).
$$
So, our inequality actually reads
$$
\frac{|v^s|^2}{2}-\partial_t((t-T)(\psi^s-\frac{|v^s|^2}{2}))
-\psi
-(t-T)v^s\cdot\nabla(\psi^s-\frac{|v^s|^2}{2})
$$
$$
-\partial_t((t-T)v^s)\cdot v^s-(t-T)\cdot\nabla v^s\cdot(v^s\otimes v^s)
-(t-T)v^s\cdot E\ge 0,
$$
which simply reduces to
$$
-\partial_t((t-T)\psi^s)
-\psi-(t-T)v^s\cdot (\nabla\psi^s+E)\ge 0.
$$
At this stage, a natural choice for $\psi$ (which has not been made so far) is
\begin{equation}\label{ansatz psi}
\psi=-\partial_t((t-T)\psi^s)
\end{equation}
so that now our inequality just reads
$$
v^s\cdot (\nabla\psi^s+E)\ge 0.
$$
So we are finally left with minimizing in $E$
$$
\int_{Q} \epsilon E\cdot \nabla\psi+\frac{\epsilon|E|^2}{2}
$$
where $E$ is subject to 
$$
v^s(t,x,a)\cdot(\nabla\psi^s+E)(t,x)\ge 0,\;\;\;\forall (t,x,a)\in Q'.
$$
According to assumption (\ref{IVP-absorbing}), which means that the velocity field $v^s$ is ``weakly absorbing",
$E=-\nabla\psi^s$ turns out to be the only possible choice for $E$!
Indeed, once $(t,x)$ is fixed, let us set $V=-(\nabla\psi+E)(t,x)$, which must satisfy
$$
v^s(t,x,a)\cdot V\le 0,\;\;\;\forall a \in\mathcal A.
$$
 By (\ref{IVP-absorbing}), we can write
$$
V=\int_{a\in \mathcal A}v^s(t,x,a)\lambda(da)
$$ 
for some nonnegative measure $\lambda$ (that, of course, depends on $t,x$ and $V$)
and, therefore, 
$$
|V|^2=
\int_{a\in \mathcal A}V\cdot v^s(t,x,a)\lambda(da)\le 0
$$
which implies $V=0$, i.e. $E(t,x)=-\nabla\psi^s(t,x)$.
So, the value of $\mathcal J$ is just
$$
\mathcal J=
\int_{Q} -\epsilon \nabla\psi^s\cdot \nabla\psi+\frac{\epsilon|\nabla\psi^s|^2}{2}.
$$
Since $\psi=-\partial_t((t-T)\psi^s)$, we have
$$
\int_{Q} -\nabla\psi^s\cdot \nabla\psi=\int_{Q} \nabla\psi^s\cdot \partial_t((t-T)\nabla\psi^s)
=\int_{Q} |\nabla\psi^s|^2+(t-T)\partial_t(\frac{|\nabla\psi^s|^2}{2})
$$
$$
=\int_{Q} \frac{|\nabla\psi^s|^2}{2}+T\int_{x\in D} \frac{|\nabla\psi^s(0,x)|^2}{2}
$$
so that
\begin{equation}\label{Jquart}
\mathcal J
=\int_{Q} \epsilon|\nabla\psi^s|^2+T\int_{x\in D} \frac{\epsilon|\nabla\psi^s(0,x)|^2}{2}
\end{equation}
After the computation of $\mathcal J$ we have just performed, we still need
an evaluation of both $BT(\theta,A)$ and $\int_{Q} \psi$
in our lower bound (\ref{lower}).
Concerning the first one, we have, by definition (\ref{IVP-BT}),
$$
BT(\theta,A)=-
\int_{(x,a)\in D\times\mathcal A} c_0(x,a)\theta(0,x,a)+q_0(x,a)\cdot A(0,x,a)
$$
$$
=T\int_{(x,a)\in D\times\mathcal A} c_0(x,a)(\psi^s(0,x)-\frac{|v^s(0,x,a)|^2}{2})+q_0(x,a)\cdot v^s(0,x,a)
$$
(since $\theta$ and $A$ are given by (\ref{ansatz}))
$$
=T\int_{(x,a)\in D\times\mathcal A} c_0(x,a)(\psi^s(0,x)+\frac{|v^s(0,x,a)|^2}{2})
$$
(since $c_0v^s(0,\cdot)=q_0$).
Because of (\ref{MEPS}), we have $1+\epsilon\Delta\psi^s=\int_a c$ and, therefore,
$$
\int_{(x,a)\in D\times\mathcal A}
 c_0(x,a)\psi^s(0,x)=\int_{x\in D} (1+\epsilon\Delta\psi^s(0,x))\psi^s(0,x)
=\int_{x\in D} \psi^s(0,x)-\epsilon|\nabla \psi^s(0,x)|^2.
$$
Thus
$$
BT(\theta,A)=T\int_{x\in D} \psi^s(0,x)-\epsilon|\nabla \psi^s(0,x)|^2+T\int_{(x,a)\in D\times\mathcal A}c_0(x,a)\frac{|v^s(0,x,a)|^2}{2}.
$$
Let us now move to the second term $\int_{Q} \psi$. From (\ref{ansatz})), we get
$$
\int_{Q} \psi=\int_Q -\partial_t((t-T)\psi^s)=-T\int_{x\in D}\psi^s(0,x).
$$
So, we have found
$$
BT(\theta,A)+\int_{Q} \psi=
-T\int_{x\in D}\epsilon|\nabla \psi^s(0,x)|^2+T\int_{(x,a)\in D\times\mathcal A}c_0(x,a)\frac{|v^s(0,x,a)|^2}{2}.
$$
Using the value of $\mathcal J$ given by (\ref{Jquart}), we see that our lower bound (\ref{lower}) now reads
\begin{equation}\label{lowerbis}
J(c_0,q_0)\ge 
T\int_{(x,a)\in D\times\mathcal A}c_0(x,a)\frac{|v^s(0,x,a)|^2}{2}.
+\int_{Q} \epsilon|\nabla\psi^s|^2-T\int_{x\in D} \frac{\epsilon|\nabla\psi^s(0,x)|^2}{2}
\end{equation}
\subsubsection*{Final step}
Collecting the results of Step 1 and Step 2, namely (\ref{upper}) and (\ref{lowerbis}), we have finally obtained
$$
0\le I(c_0,q_0)- J(c_0,q_0)\le
-T\int_{(x,a)\in D\times\mathcal A}c_0(x,a)\frac{|v^s(0,x,a)|^2}{2}
+T\int_{x\in D} \frac{\epsilon|\nabla\psi^s(0,x)|^2}{2}
$$
$$
-\int_Q\frac{\epsilon|\nabla\psi^s|^2}{2}+\int_{Q'} \frac{c|v^s|^2}{2}.
$$
As a matter of fact, the right-hand side is exactly zero because of the conservation of energy (\ref{MEPS-energy}).
So there is no duality gap and the proof of Theorem \ref{nogap} is now complete.

\subsection*{Acknowledgment}
This work originated from a visit of the first author
at the Department of Pure Mathematics and Mathematical Statistics (DPMMS),
University of Cambridge, in May 2018, supported by the grant MAFRAN.
A part of the work was performed while he was visiting the Hausdorff Research Institute for Mathematics, Bonn,
in 2019, during the program: ``Interfaces and Instabilities in Fluid Dynamics" co-organized by L. Sz\'ekelyhidi Jr. and G. Weiss.

\end{document}